\documentclass[preprint,12pt]{elsarticle}

 \usepackage{graphics}
 
\usepackage{booktabs} 
\usepackage{array} 
\usepackage{paralist} 
\usepackage{verbatim} 
\usepackage{subfigure} 
\usepackage{pstricks-add}
\usepackage{amsmath}
\usepackage{amsfonts}
\usepackage{amssymb}
\usepackage{graphicx} 

\usepackage[left]{lineno} 

\usepackage{multicol}

\usepackage{cmbright} 

\usepackage{geometry} 

\usepackage{color}
\usepackage{parskip}

\usepackage{amssymb}

\def\qed{\ifmmode\square\else\nolinebreak\hfill
$\square$\fi\par\vskip12pt}

\newtheorem{lemma}{Lemma}
\newtheorem{obs}{Observation}

\newtheorem{theorem}{Theorem}
\newtheorem{definition}{Definition}
\newtheorem{example}{Example}
\newtheorem{claim}{Claim}

\journal{Discrete Applied Mathematics}

\begin{document}

\begin{frontmatter}



\title{The chromatic number of comparability 3--hypergraphs}


\author{Natalia  Garc\'{i}a-Col\'{i}n
\\
{\normalsize Instituto de Matem\'{a}ticas, UNAM}\\
{\normalsize natalia.garciacolin@im.unam.mx} 
\\[+12pt] 
Amanda Montejano
\\
{\normalsize UMDI Facultad de Ciencias, UNAM}\\
{\normalsize amandamontejano@ciencias.unam.mx} 
\\[+12pt] 
Deborah Oliveros 
\\
{\normalsize Instituto de Matem\'{a}ticas, UNAM}\\
{\normalsize dolivero@matem.unam.mx}}

\begin{abstract}
Beginning with the concepts of orientation for a $3$--hypergraph and transitivity for an oriented $3$--hypergraph, it is natural to study the class of  comparability $3$--hypergraphs (those that can be  transitively oriented). In this work we show three different behaviors in respect  to the relationship between the chromatic number and the clique number of a comparability $3$--hypergraph, this is in contrast with the fact that a comparability simple graph is a perfect graph.

\end{abstract}


\end{frontmatter}


\section{Introduction and motivation}

In \cite{Cyclic2013} the authors introduce the concepts of orientation for a $3$--hypergraph, transitivity for an oriented $3$--hypergraph,  and define the class of \emph{comparability $3$--hypergraphs} as the class of non oriented $3$--hypergraphs, which can be transitively oriented (precise definitions are provided in Section \ref{sec:pre}). These  $3$--hypergraphs are a natural generalization of (simple) comparability graphs (graphs which can be transitively oriented or, equivalently,  graphs associated to a partially ordered set).

Comparability graphs are well known to be perfect graphs. A graph is said to be \emph{perfect} if all of its induced subgraphs have chromatic number equal to their clique number. This concept was introduced by Claude Berge in 1961 \cite{Berge1961}. In essence it means  that a graph is perfect if for every of its induced subgraphs the chromatic number is  as low as possible  in terms of its clique number. Thus, it is natural to ask whether or not comparability $3$--hypergraphs are perfect in this sense. 

Hypergraphs have been studied in relation to  perfection in \cite{Lovasz1972253}, 
and \cite{lasvergnas}. However,  the precise concept of  perfection for hypergraphs  remains imprecise, to the best of our knowledge. Aiming to find a suitable definition of perfection in hypergraphs, we study the relationship between the chromatic number  and the clique number  of comparability $3$--hypergraphs. 

We define  a \emph{$3$--hypergraph} $H$ as usual,  $H=(V(H),E(H))$ where $V(H)$ is the set of \emph{vertices} of $H$, and $E(H)\subseteq$$ {V(H)}\choose{3}$ is the set of \emph{edges}. The  \emph{chromatic number}, $\chi(H)$, is defined as the minimum $k$, such that $V(H)$ can be partitioned into $k$ parts, called \emph{color classes}, in such a way that no edge of $H$ is monochromatic, in other words, no edge is contained in a single color class. The   \emph{clique number}, $\omega (H)$, of a $3$--hypergraph $H$ is the largest cardinality of a subset of $V(H)$ inducing a complete $3$--hypergraph. 

Given that for any complete $3$--hypergraph on $n$ vertices,$K^3_n$, we have that $\chi(K^3_n)=\left\lceil \frac{n}{2}  \right\rceil $, then for any $3$--hypergraph the following equation holds:
\begin{equation}
\left\lceil \frac{\omega(H)}{2}  \right\rceil \leq \chi(H).\end{equation}

In this paper we study three important subclasses of comparability $3$--hypergraphs which show three different behaviors in relation to (1).

Firstly, we exhibit a family of  comparability $3$--hypergraphs for which the difference,  $\chi(H) - \left\lceil \frac{\omega(H)}{2}  \right\rceil$, is arbitrarily large. 

Secondly, we exhibit an interesting   subclass of comparability $3$--hypergraphs, named cyclic permutation $3$--hypergraphs (the analogues of permutation graphs), such that their chromatic number is bounded by a (linear) function of its clique number. 

Finally, we exhibit another interesting subclass of comparability $3$--hypergraphs namely, the ones associated to a family of intervals in the circle. For these hypergraphs the chromatic number is  as low as it can be in respect to their clique, that is, equality holds in (1).

The paper is organized as follows: in Section \ref{sec:pre} we state the required definitions and preliminary results necessary to prove our main theorems, our main results are stated in Section \ref{sec:thm} and the proofs are located in the remaining sections.


\section{Definitions and preliminaries}\label{sec:pre}

Let $X$ be any set of order $n$. A \emph{linear ordering} of $X$ is a bijection $\phi :\{1,2,...,n\} \to X$. A \emph{cyclic ordering} of $X$ is an equivalent class of the set of linear orderings  with respect to the \emph{cyclic equivalence relation}  defined as: $\phi\sim \psi$, if and only if there exists $k\leq n$, such that $\phi(i)= \psi(i+k)$ for every $i\in \{1,2,...,n\}$ where $i+k$ is taken modulo $n$. For the remainder of this paper we will denote each cyclic ordering, $[\phi]$, in cyclic permutation notation, $\big(\phi(1)\, \phi(2)\, \ldots\,\phi(n)\big)$. For example, there are two different cyclic orderings of $\{u,v,w\}$, namely $(u\ v\ w)$ and $(u\ w\ v)$, where $(u\ v\ w)=(v\ w\ u)=(w\ u\ v)$ and $(u\ w\ v)=(v\ u\ w)=(w\ v\ u)$.

Given a $3$--hypergraph $H$, an \emph{orientation} of $H$ is an assignment of exactly one of the two possible cyclic orderings  to each of its edges. An orientation of a $3$--hypergraph is called an \emph{oriented $3$--hypergraph}, and we denote the oriented edges by $O(H)$. 

\begin{figure}[t]
\begin{center}
\scalebox{0.15}{\includegraphics{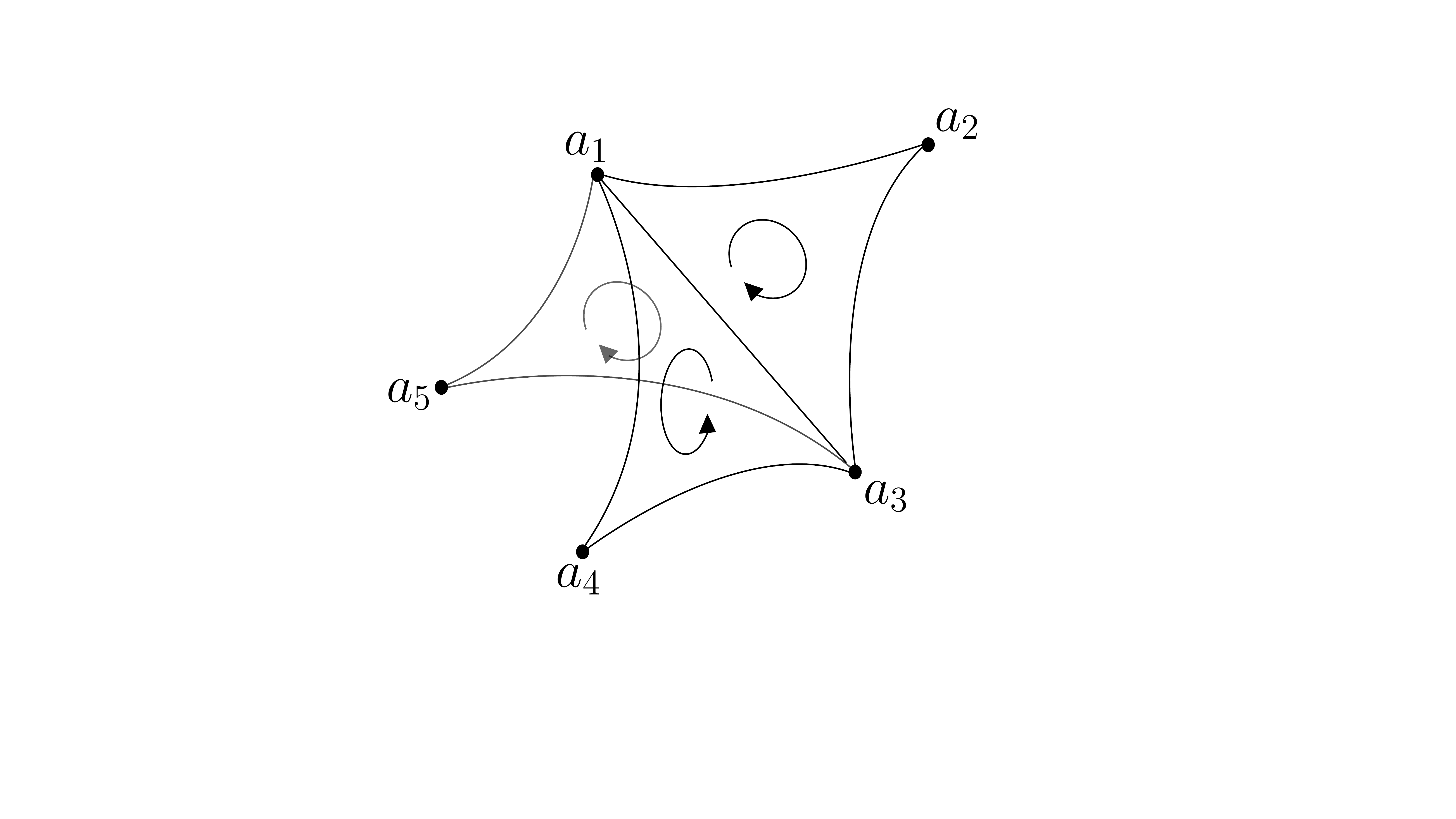}} \\[-.80in]
\vspace{.5cm}
\caption{A non transitive oriented $3$--hypergraph for which its underlying graph is a  comparability $3$-hypergraph.} \label{fig1}
\end{center}
\end{figure}

\begin{example}\label{ex}
Let $H=(V(H),E(H))$ be a $3$--hypergraph  with $V(H)=\{ a_1,a_2, a_3, a_4, a_5\}$ and 
$E(H)= \{ \{a_1,a_2, a_3\} ,\{a_1, a_3, a_4\} , \{ a_1,a_3, a_5\} \}$,
then a possible orientation of $H$ could be $O(H)=\{ (a_1\,a_2\,a_3), (a_1\,a_4\,a_3),(a_1\,a_3\,a_5) \} $ obtaining the oriented $3$--hypergraph depicted in  Figure \ref{fig1}.
\end{example}

It is usual to associate transitive oriented graphs to partial (linear) orders.  Similarly, we can associate to partial cyclic orders (a ternary relation which is cyclic, asymmetric and transitive)  transitive oriented $3$--hypergraphs in the following manner:

\begin{definition}
An oriented $3$--hypergraph $H$ is said to be \emph{transitive}, if whenever $(u\, v\, z)$ and $(z\, v\, w) \in O(H)$ then $(u\, v\, w) \in O(H)$ (this implies also $(u\, w\, z) \in O(H)$). 
\end{definition}

Now it is natural to define and study the following class of  $3$-hypergraphs. 

\begin{definition}\label{comphyp}
A non-oriented $3$--hypergraph  is called a \emph{comparability $3$--hypergraph} if it admits a transitive orientation. 
\footnote{In  \cite{Cyclic2013} the authors defined this class as ``cyclic comparability $3$-hypergraphs" however we believe that it should be simply called comparability $3$-hypergraphs according to the classical concept of comparability graphs. }
\end{definition}

The oriented $3$--hypergraph defined in Example \ref{ex} is not transitive, however, its underlaying $3$--hypergraph $H$ is a comparability $3$--hypergraph since it can be  transitively oriented;  take for instance $O'(H)=\{(a_1\,a_3\,a_2), (a_1\,a_3\,a_4),(a_1\,a_3\,a_5) \}$. In contrast, a $3$--hypergraph with four vertices and three edges is not a comparability $3$--hypergraph.

An oriented $3$--hypergraph containing all possible $3$--edges is called a \emph{$3$--hypertournament}.  As in the case of oriented graphs, there is a unique transitive oriented $3$--hypertournament with $n$ vertices that we will denote by $TT_n^3$  \cite{Cyclic2013}.

In this paper we study the chromatic number of two important subclasses of  comparability $3$--hypergraphs, namely the ones associated to a cyclic permutation, and the ones associated to a finite family of closed intervals in the circle. We will now proceed to define both classes.

\subsection{Cyclic Permutation $3$--hypergraphs}

A \emph{cyclic permutation} is a cyclic ordering of the set $\{1,2,...,n\}$. This is, an equivalence class $[\phi]$ of the set of bijections $\phi:\{1,2,...,n\} \rightarrow \{1,2,...,n\}$, in respect to the cyclic equivalence relation.  Let $[\phi]$ be a cyclic permutation. Three elements $i, j, k \in \{1,2,...,n\}$, with $i< j< k$, are said to be in \emph{clockwise order} in respect to $[\phi]$, if there is $\psi\in [\phi]$, such that $\psi^{-1}(i)<\psi^{-1}(j)<\psi^{-1}(k)$. Otherwise the elements $i,j,k$ are said to be in \emph{counter-clockwise order} with respect to $[\phi]$. 

\begin{definition}
The \emph{oriented $3$--hypergraph $H_{[\phi ]}$ associated to a cyclic permutation} $[\phi ]$ is the $3$--hypergraph with vertex set $V(H_{[\phi ]})=\{1,2,...,n\}$ where the edges are the triplets $\{i,j,k\}$ with $i<j<k$, which are in clockwise order in respect to $[\phi ]$, and where the edge orientations are induced by $[\phi ]$.
\end{definition}

 It is not difficult to check that a $3$--hypergraphs associated to a cyclic permutations is indeed transitive, and thus its underlying $3$--hypergraph is a comparability $3$--hypergraph.

\begin{example}
Consider  the identity cyclic permutation $(1\,2\,\dots\, n)$ and its reversed cyclic permutation $(n\,\dots \,2 \,1)$, then their respective associated oriented $3$--hypergraphs are the transitive $3$--hypertournament $TT_n^3$ and  the null $3$--hypergraph on $n$ vertices respectively. 
\end{example}

A non-oriented $3$--hypergraph is  called a \emph{cyclic permutation $3$--hypergraph}, if it can be oriented in such a way that the resulting oriented hypergraph is isomorphic to the one associated to a cyclic permutation. Emulating a classical result of Pnueli, Lempel and Even \cite{ple}, the authors of  \cite{Cyclic2013} provided a characterization of the class of cyclic permutation $3$--hypergraphs  in terms of comparability $3$--hypergraphs. 

\begin{theorem}[Garcia-Colin, Montejano, Montejano, Oliveros \cite{Cyclic2013}] \label{thm:gmmo}
$H$ is a cyclic permutation $3$--hypergraph, if and only if $H$  and its complement $\overline{H}$ are  comparability $3$--hypergraphs.
\end{theorem}

\subsection{Circle interval $3$--hypergraphs} 

A simple graph $G$ is called an \emph{interval graph}, if there exists a finite set of closed intervals in the real line, whose intersection graph is $G$.  Interval graphs,  as well as its complements are well known to be perfect graphs.  In a similar  way, we will study $3$--hypergraphs arising  from  finite sets of closed intervals in the circle $S^1$.

\begin{definition}
The \emph{$3$--hypergraph $H_{\mathcal{F}}$ associated to a finite family of closed intervals in the circle $\mathcal{F}$}, is the $3$--hypergraph with vertex set $V(H_{\mathcal{F}})=\mathcal{F}$ whose edges are the triplets of vertices with the property that their corresponding intervals are pairwise disjoint. 
\end{definition}

For instance, if $\mathcal{F}$ is a family of $n$ pairwise disjoint closed intervals in the circle, then $H_{\mathcal{F}}$ is the complete $3$--hypergraph on $n$ vertices.

A $3$--hypergraph $H$ is called a \emph{circle interval  $3$--hypergraph}, if there exists a set of  closed intervals in the circle for which the associated $3$--hypergraph is isomorphic to $H$. It is not difficult  to prove that any circle interval  $3$--hypergraph $H$ is a comparability $3$--hypergraph, since the natural orientation of $H$,  given by the clockwise order in which each of the intervals  appears in $S^1$, is transitive.

\section{Main results}\label{sec:thm}

Recall that comparability simple graphs are always perfect. We will show  three classes of  comparability $3$--hypergraphs exhibiting three different types of behaviors among the relationship between $\chi$ and $\omega$, thus proving that comparability $3$--hypergraphs are far from perfection.

First we will show that within the class of  comparability $3$--hypergraphs the difference  $\chi-\left\lceil \frac{\omega}{2}  \right\rceil$ is not bounded.

\begin{theorem} \label{thm:main1}
For each positive integers $w$ and $k$, such that $\left\lceil \frac{w}{2}  \right\rceil \leq k$,
 there exists a  comparability $3$--hypergraph with clique number $w$ and chromatic number at least 
 $k$.
\end{theorem}

In view of the above result it is an interesting problem to study families of comparability $3$--hypergraphs for which the chromatic number of each member in the family is bounded by a function of its clique number.

\begin{theorem} \label{thm:main2}
Let $H$ be a cyclic permutation $3$--hypergraph,  then $\chi(H) \leq \omega(H)-1$. Furthermore, this bound is tight.
\end{theorem}

Finally, we show an infinite family of comparability $3$--hypergraphs for which the the chromatic number is as small as it can be in terms of the clique number.

\begin{theorem} \label{thm:main3}
Let $H$ be a circle interval $3$-hypergraph, then $\chi(H)=\left\lceil{ \frac {\omega(H) }{2}}\right\rceil$.
\end{theorem}

\section{Comparability $3$--hypergraphs with large chromatic number and fixed clique number}

In this section we will prove Theorem \ref{thm:main1} by constructing a comparability $3$--hypergraph with fixed clique number  and   arbitrarily large chromatic number. We will use a construction given in  \cite{Mycielskitype2013} and prove that the $3$--hypergraphs obtained can be transitively oriented.

Let $H$ be a $3$--hypergraph of order $n$, with vertex set $V(H)$ and edge set $E(H)$. We define $\mu (H)=(V,E)$ as follows:

\begin{itemize}
\item[i)] $V=V_1\cup V_2 \cup \{w\}$ where $V_1:= \{(v,1)|v\in V(H)\}$ and $V_2:= \{(v,2)|v\in V(H)\}$ 
\item[ii)] $E=E_1\cup E_2$ where:

$E_1$ is the set of triplets in $V_1\cup V_2$ whose projection (on to H) is an edge of $H$, that is,  \\
$\{(u_{1},t_1),(u_{2},t_2),(u_{3},t_3)\}\in E_1$,  if and only if
$\{u_{1},u_2 ,u_{3}\}\in E(H)$ and $t_i\in\{1,2\}$.

$E_2$ is formed by the pairs of $V_1\cup V_2$ whose projection (on to H) is a
single vertex, together  with the vertex $w$, that is, 
$\{(u_1 , 1),(u_2, 2),w \}\in E_2$  if and only if $u_1=u_2$.

 \end{itemize}
 
Note that $|V(\mu (H))|=|V|=2n+1$, $|E_2|=n$, $|E_1|=8|E(H)|$,  and that the subhypergraphs induced by $V_1$ and $V_2$ are both copies of $H$ Furthermore:
\begin{equation}
\omega (\mu(H))=\omega(H) \  { and } \ \chi (\mu (H))=\chi (H)+1
\end{equation}

See \cite{Mycielskitype2013} for a complete proof of this fact. 

Next we will show that, if $H$ is a comparability $3$--hypergraph then  $\mu (H)$ is also a   comparability  $3$--hypergraph. Let $O(H)$ be a transitive orientation of $H$. We will define a natural transitive orientation $O(\mu(H))$ using the following rules:

Let $e=\{ x,y,z\}\in E(\mu(H))$. 
\begin{enumerate}[i)]
\item If $e\in E_2$, say $x=(u_{1},1)$, $y=(u_{1},2)$, and $z=w$, then
orient $e$ simply as $(x\ y\ z)$.  
\item If  $e \in E_1$, say $x=(u_{1},t_1)$, $y=(u_{2},t_2)$, $z=(u_{3},t_3)$ where $\{u_{1},\ u_2 \ ,u_{3}\}\in E(H)$, and $t_i\in\{1,2\}$, then orient  $\{ x,y,z\}$ as induced by the orientation of $\{u_1\ u_2\  u_3\}$, i.e.  $(x \ y \ z) \in O(\mu(H)) $ if and only if $(u_1\ u_2\  u_3)\in O(H)$. 
\end{enumerate}

\begin{claim}
$O(\mu(H))$  is transitive.  
\end{claim}
 
{\it{Proof}}
[Proof of Claim 1] We have to check that the definition of transitivity is satisfied by every pair of edges. Notice first, that if $e$ and $f$ are two edges in $\mu(H)$,  and one of them is in $E_2$ then $|e \cap f |=1$, so there is no conflict between the orientations of $e$ and $f$. 

Suppose $e, f$ are such that $|e \cap f |=2$ then, by the previous observation, both $e$ and $f $ are in $E_1$. First suppose $e=(u\ v \ z )$,  $f=(v\ z\  x)\in O(\mu(H))$, then the orientation satisfies the transitivity condition.Secondly, suppose $e=(u\ v \ z )$,  $f=(z\ v\  x)\in O(\mu(H))$,  $u=(u_1,t_1)$, $v=(v_1,t_2)$, $z=(z_1,t_3)$ and $x=(x_1,t_4)$ where $t_i\in \{ 1,2\}$ for $i\in \{ 1,2,3,4\}$.
By the orientation rules $(u_1\  v_1\  z_1)$ and $(z_1\ v_1\  x_1)\in O(H)$. By the definition of transitivity of $H$ it follows that 
$(u_1\  v_1\  x_1)\in O(H)$,  implying that $(u\ v\  x)\in O(\mu(H))$, thus the definition of transitivity is satisfied.
\qed

Claim 1 proves that $\mu(H)$ is a  comparability $3$--hypergraph.\\


\textbf{Proof of Theorem \ref{thm:main1}.}
Start with $H$ a comparability $3$--hypergraph with clique number $w$ and chromatic number $\left\lceil \frac{w}{2}  \right\rceil$ (take for instance the complete $3$--hypergraph with $w$ vertices). Define $\mu^{k+1}(H)=\mu(\mu^k(H))$ where $\mu^1(H)=\mu(H)$. Then, using the facts discussed in the previous paragraphs, for any $k$, $\mu^k(H)$ is a comparability $3$--hypergraph with clique number $w$ and chromatic number at least $\left\lceil \frac{w}{2}  \right\rceil+k$. \qed

\section{Even $3$--hypergraphs}

In this section we prove the upper bound  of Theorem \ref{thm:main2}  using  an interesting class of $3$--hypergraphs called even $3$--hypergraphs. 

\begin{definition} A $3$--hypergraph $H$ is said to be \emph{even}, if every four vertex subset of $V(H)$ induces  an even number of edges.
\end{definition}

The above definition arises naturally when studying the class of $3$--hypergraphs associated to cyclic permutations. In fact, the class of cyclic permutation $3$--hypergraphs is  an important subclass of the class of even $3$--hypergraphs. To see this recall that by  Theorem \ref{thm:gmmo} a cyclic permutation $3$--hypergraph $H$ is such that both $H$ and its complement $\overline{H}$ can be transitively oriented. Consequently, any subset  of four vertices in $V(H)$ induces $0$, $1$, $2$ or $4$ edges in $H$ as well as in $\overline{H}$. Therefore,

\begin{obs}\label{rem:even}
Any cyclic permutation $3$--hypergraph is an even $3$--hypergraph.
\end{obs}

Next we prove two lemmas which relate the chromatic number (respectively, the clique number) of an even $3$--hypergraph $H$ with the chromatic number  (respectively, the clique number) of a simple graph associated to $H$.

Let $H$ be a $3$--hypergraph and $v \in V(H)$; we define the \emph{link} of $v$ as the simple graph $link_H(v)$, where $ V(link_H(v))= \{ w \neq v \in V(H) | \{v, w\} \subset e, \text{ for some } e \in E(H) \}$, and $ E(link_H(v))= \{ \{u,w\} \subset V(link_H(v)) | \{v,u, w\} \in E(H) \}$.

\begin{lemma} \label{lem:link_chi} Let $H$ be an even $3$--hypergraph, then $\chi(H)\leq \chi(link_H(v))$ for all $v \in V(H)$.
\end{lemma}

{\it{Proof}} For a fixed  $v \in V(H)$, let  $\chi(link_H(v))=k$ and $c: V(link_H(v)) \rightarrow \{1,...,k\}$ be a proper coloring. We  extend this coloring to a coloring of $H$ by defining $c': V(H) \rightarrow \{1,...,k\}$ as $c'(u)=c(u)$ for $ u \in V(link_H(v))$, and $c'(u)=k$ for $u \in V(H) \setminus V(link_H(v))$. In order to prove that $\chi(H)\leq k$ it remains to prove that $c'$ is a proper coloring of $H$, that is,  there are no edges of $H$ whose vertices have the same color.

Let $V(H)=V_1\cup V_2\cup...\cup V_k$ be the partition induced by the coloring defined in the previous paragrph.  Suppose, by contradiction, that $\{u_1,u_2,u_3\}\in E(H)$ is such that $\{u_1,u_2,u_3\}\subseteq V_i$ for some $i\in\{1,...,k\}$. Firstly assume that $v\in \{u_1,u_2,u_3\}$, then the pair $\{u_1,u_2,u_3\}\setminus v$ is an edge of $link_H(v)$, which is a contradiction since the coloring restricted to $V(link_H(v))$ is proper. Assume now that $v\not\in \{u_1,u_2,u_3\}$, and consider the set of vertices $\{v,u_1,u_2,u_3\}$, which  should induce an even number of edges of $H$. Thus, at least one pair in $\{u_1,u_2,u_3\}$ is an edge of $link_H(v)$,  which is a contradiction for the same reason.  Hence,  $c'$ is a proper $k$--coloring of $H$ and the statement follows true.
\qed

\begin{lemma}\label{lem:link_omega}
Let $H$ be an even $3$--hypergraph then $\omega(H) \geq \omega(link_H(v))+1$ for all $v \in V(H)$.
\end{lemma}

{\it{proof}}
Set $\omega(link_H(v))=\omega$, and let  $K_{\omega}$ be a clique of $link_H(v)$,  then the subhypergraph of $H$ induced by the vertex set $\{v\} \cup V(K_{\omega})$ is a clique of $H$ due to the eveness of $H$.
\qed

Now we are ready to prove  the upper bound  of Theorem \ref{thm:main2}.

\textbf{Proof of Theorem \ref{thm:main2} (upper bound)}:
Let $H$ be a cyclic permutation $3$--hypergraph. By Observation \ref{rem:even}, $H$ is an even $3$--hypergraph. Hence, by Lemma \ref{lem:link_chi} and Lemma \ref{lem:link_omega} we obtain $\chi(H)\leq \chi(link_H(v))$ and $\omega(link_H(v))\leq \omega(H) - 1$  respectively. To conclude the proof it remains to argue that the graph $link_H(v)$ is a perfect graph, thus $\chi(link_H(v))=\omega(link_H(v))$ and the statement holds true.

To see that $link_H(v)$ is a perfect graph we will prove that it is a comparability graph (a graph that can be transitively oriented). We claim that a transitive orientation of $H$ naturally induces a transitive orientation of $link_H(v)$. Let $O(H)$ be a transitive orientation of $H$. Now consider the orientation of $link_H(v)$ defined as follows: for $\{u, w\} \in E(link_H(v))$ let $(u,w) \in A(link_H(v))$,  if $(v\, u\, w) \in O(H)$. By the definition of transitivity for $H$ this orientation of $link_H(v)$ is transitive. This is,  if $(u_1,u_2)$ and $(u_2,u_3)$ are arcs of $link_H(v)$ then $(u_3\, v\, u_2)$ and $(u_2\, v\, u_1)\in O(H)$, therefore $(u_3\, v\, u_1)\in O(H)$ and so $(u_1,u_3)\in A(link_H(v))$.  \qed

\section{Winding permutation $3$--hypergraphs}

In this section we prove the tightness of Theorem \ref{thm:main2} by defining a special family of cyclic permutation $3$--hypergraphs which we call winding permutation $3$-hypergraphs. 
Before proceeding we provide some terminology. For  $S=\{s_1,...s_m\}\subseteq \{1,...,n\}$ we denote as  $[\phi_S]$ the \emph{cyclic sub-permutation} of  $[\phi]$ induced by $S$, that is, $\big(\phi_S(1)\,\,...\,\phi_S(m)\big)=\big(\phi(s_1)\,\,...\,\phi(s_m)\big)$. Note that $[\phi_S]$ is a an equivalent class of the set of mappings  (not necesarly surjective) $\phi_S:\{1,\,...\, m\}\to \{1,\,...\, n\}$  in respect to the cyclic equivalence relation. If a cyclic sub-permutation  $[\phi_S]$ is such that  there is $\psi\in [\phi_S]$ with $\psi(1)<\psi(2)\,...<\psi(m)$ then we say that $[\phi_S]$ is a
\emph{clockwise increasing} cyclic sub-permutation of $[\phi]$, if there is $\psi\in[\phi_S]$ with $\psi(1)>\psi(2)\,...>\psi(m)$ we say that  $[\phi_S]$  is a \emph{clockwise decreasing} cyclic sub-permutation of $[\phi]$. With the aid of this definition we will remark that:

\begin{obs}\label{obs:omwgachi}
Let $[\phi ]$ be a cyclic permutation then:

\begin{enumerate}[i]
\item[i)] $\omega(H_{[\phi ]})$ is the length of the longest clockwise increasing cyclic sub-permutation of $[\phi]$.

\item[ii)] $\chi(H_{[\phi ]})$ is  the minimum number of clockwise decreasing cyclic sub-permutations of $[\phi]$, which cover $[\phi]$, that is, the minimum number k such that there is a partition $[n]= S_1 \cup S_2\, \ldots \cup S_k  $ which satisfies that each $[\phi_{S_i}]$ is a clockwise decreasing cyclic sub-permutation.
\end{enumerate}

\end{obs}

For example, consider the cyclic permutation $3$--hypergraph $H_{[\phi]}$ where $[\phi]=(5\,2\,6\,3\,7\,4\,1)$. Then, $\omega (H_{[\phi]})=4$ since  $(5\,6\,7\,1)$ is a clockwise increasing cyclic sub-permutation of $[\phi]$ and there are no clockwise increasing cyclic sub-permutation of $[\phi]$ with length $5$. Also $\chi(H_{[\phi]})=3$ since the  clockwise decreasing cyclic sub-permutations $(5\,2)$, $(6\,3)$ and $(7\,4\,1)$ cover $[\phi]$, and this can not be done with just two clockwise decreasing cyclic sub-permutations of $[\phi]$. 

\begin{definition} For positive integers  $q$ and $r$ we define \emph{the winding permutation $3$--hypergraph}, $W_{q, r}$, as the $3$--hypergraph with $n=r(q-1)+1$ vertices associated to the cyclic permutation $[\phi _{q, r}]$ defined as follows. For $i\in \{1,2,...,n\}$ let:

\[\phi_{q, r}(i) = \left\{ \begin{array}{l|l}
        1 & \mbox{if \hspace{.3cm}$i=r(q-1)+1=n$}\\
          2+(q-1)(r-j)+ \frac{i-j}{r}& \mbox{if \hspace{.3cm}$i\equiv j$ (mod $r$)  \hspace{.2cm}where  \hspace{.2cm}$1\leq j\leq r$}
        \end{array} \right. \]

\end{definition}

For example, let $q=5$ and $r=3$ then we have $n=13$ vertices, and $W_{ 5, 3}$ is the $3$--hypergraph associated to the permutation $(10\,\, 6\,\, 2\,\, 11\,\, 7\,\, 3\,\, 12\,\, 8\,\, 4\,\, 13\,\, 9\,\, 5\,\, 1 ) $.

Given $q$ and $r$ we distinguish two sets of useful cyclic sub-permutations of $[\phi _{q, r}]$, namely those induced by the sets $A_i=\{a:r(i-1)+1\leq a\leq ir\}$ where $1\leq i\leq q-1$, and those induced by the sets $B_i=\{a:a\equiv i \mbox{ (mod $r$)}\}$ where $1\leq i\leq r$. For instance, in the previous example we have: 

\begin{center}
$[A_1]=(10\,\, 6\,\, 2)$,  $[A_2]=(11\,\, 7\,\, 3)$, $[A_3]=(12\,\, 8\,\, 4)$,   $[A_4]=(13\,\, 9\,\, 5)$

$[B_1]=(10\,\, 11\,\,12 \,\,13 \,\,1)$, $[B_2]=(6\,\, 7\,\,8 \,\,9)$, $[B_3]=(2\,\, 3\,\,4 \,\,5)$ 
\end{center}

By definition, it follows that:

\begin{obs}\label{obs:ab}
Every cycle sub-permutation of the form $[A_i]$ is a clockwise decreasing  cycle sub-permutation, while every  cycle sub-permutation of the form $[B_i]$ is a clockwise increasing  cycle sub-permutation.
\end{obs}

In the next two lemmas we compute the clique and chromatic numbers of $W_{q,r}$.

\begin{lemma} \label{lem:clique} For any integers $q$ and $r$, the clique number  $\omega(W_{q, r})=q$.
\end{lemma}

{\it{Proof}}  
Note that $(1\,\,2\,...q)$ is an increasing  clockwise cyclic sub-permutation of $[\phi_{q, r}]$ thus, by Observation \ref{obs:omwgachi}, $\omega (W_{q, r}) \geq q$. 

In order to prove that $\omega (W_{q, r})\leq q$, let $(s_1\, ...\, s_{\omega})$ be a maximal clockwise increasing cyclic sub-permutation of $[\phi_{q, r}]$. First note that, by Observation \ref{obs:ab}, it follows that $|\{s_1,..., s_{\omega}\}\cap A_i|\leq 2$ for all $ i\in\{1,..., q-1\}$. Moreover, if $|\{s_1,..., s_{\omega}\}\cap A_j|= 2$ for some $j\in\{1,..., q-1\}$ then $|\{s_1,..., s_{\omega}\}\cap A_i|\leq 1$ for all $ i\in\{1,..., q-1\}\setminus \{j\}$ and $1\not\in \{s_1,...,s_{\omega}\}$. On the other hand, if $s_i= 1$ for some $1\leq i\leq \omega$, then $|\{s_1,..., s_{\omega}\}\cap A_i|\leq 1$ for every $i\in \{1,...q-1\}$. In both cases it we must have  $\omega \leq q$ completing the proof.
\qed

\begin{lemma}\label{lem:chi}
For any integers $q$ and $r$, the chromatic number  $\chi(W_{q, r})=\left\lceil \frac{n}{r+1}\right\rceil$.

\end{lemma}

{\it{Proof}}
In this proof we will denote $[\phi]=[\phi_{q, r}]$, for convenience. In order to show that $\chi (W_{q, r})\leq \left\lceil \frac{n}{r+1}\right\rceil$ we will exhibit a set of $\left\lceil \frac{n}{r+1}\right\rceil$ clockwise decreasing cyclic sub-permutations which cover $[\phi]$. 

Take the cyclic sub-permutation $[C_1]$ induced by the fist $r+1$ elements:  $(\phi(1)\,\, ... \,\, \phi(r+1))$, then take the cyclic sub-permutation $[C_2]$ induced by the next $r+1$ elements: $(\phi(r+2)\,\, ... \,\, \phi(2r+2))$,  and continue until the last cycle sub-permutation $[C_k]$ where $k=\left\lceil \frac{n}{r+1}\right\rceil$. Note that $|C_i|=r+1$ for $1\leq i\leq k-1$ and $|C_k|=m$ where $n=(k-1)(r+1)+m$. For instance, in the previous example with $q=5$ and $r=3$ we have $k=\left\lceil \frac{n}{r+1}\right\rceil=\left\lceil \frac{13}{4}\right\rceil=4$, $m=1$, and the obtained partition is: 

\begin{center}
$[C_1]=(10\,\, 6\,\, 2\,\, 11)$, $[C_2]=(7\,\, 3\,\, 12\,\, 8)$, $[C_3]=(4\,\, 13\,\, 9\,\, 5)$ and $[C_4]=(1)$
\end{center}

We claim that for all $1\leq j\leq k$ the  cycle sub-permutation $[C_j]$ is clockwise decreasing. Let $j\in\{1,...,k-1\}$ be fixed. We first note that $C_j=\{\phi\big((j-1)(r+1)+1\big)\,\, ... \,\, \phi\big(j(r+1)\big)\}$, intersects at most two consecutive sets of $A_1,A_2,...A_{q-1}$. Therefore, by Observation \ref{obs:ab}, $[C_j]$ is made up of two clockwise decreasing subsequences. Hence, it remains to show that  $\phi\big(j(r+1)\big)> \phi\big((j-1)(r+1)+1\big)$. To see this, note that both $j(r+1)$ and $(j-1)(r+1)+1$ are congruent to $j$ modulo $r$. Thus, by using the definition of  $\phi$, we obtain $\phi\big(j(r+1)\big)=2+(q-1)(r-j)+j$, and  $\phi\big((j-1)(r+1)+1\big)=2+(q-1)(r-j)+(j-1)$, being $\phi\big((j-1)(r+1)+1\big)+1=\phi\big(j(r+1)\big)$ as desired.

Recall that we set $k=\left\lceil \frac{n}{r+1}\right\rceil$. In order to prove that $\chi (W_{q, r})\geq k$ we will proceed by contradiction. Suppose that $\{S_i\}_{i=1}^{k-1}$ is a partition of \{1,2,...,n\} inducing a set of clockwise decreasing cyclic sub-permutations $[S_1]$, $[S_2]$,..., $[S_{k-1}]$, which covers $[\phi]$. Then, by the pigeonhole principle, for some  $j\in\{1,2,...k-1\}$ it must be true that $|S_j|\geq r+2$. Since, by Observation \ref{obs:ab}, it follows that $|S_j\cap B_i|\leq 2$ for all $ i\in\{1,..., r\}$, and moreover, if $|S_j\cap B_i|= 2$ for some $i\in\{1,..., r\}$ then $|S_j\cap B_i|\leq 1$ for all $ i\in\{1,..., r\}\setminus \{i\}$, which is a contradiction. This completes the proof.
\qed

\textbf{Proof of Theorem \ref{thm:main2} (tightness)}: For any  given $q\geq3$ we will exhibit a cyclic permutation $3$--hypergraph $H$, with $\omega(H)=q$ and $\chi(H)=q-1$. Take $H=W_{q,r}$
with $r>q-3$. By Lemma \ref{lem:clique} and Lemma \ref{lem:chi} we know that  $\omega(H)=q$ and $\chi(H)=\left\lceil \frac{n}{r+1}\right\rceil$ where $n$ is the order of $H$ and satisfies $n=r(q-1)+1$. Consequently, $\chi(H)=\left\lceil \frac{r(q-1)+1}{r+1}\right\rceil$, which can be written as $\chi(H)=\left\lceil (q-1)-\frac{q-2}{r+1}\right\rceil$. It remains to observe that $0<\frac{q-2}{r+1} <1$ by hypothesis, then $\chi(H)=q-1$ as claimed. 
\qed

\section{Circle interval $3$--hypergraphs}

Recall that a circle interval $3$-hypergraph $H_{\mathcal{F}}$  associated to a  finite family ${\mathcal{F}}$ of closed intervals in the circle, is a $3$--hypergraph with vertex set $V(H_{\mathcal{F}})=\mathcal{F}$ and  whose edges are  the triplets of vertices whose corresponding intervals are pairwise disjoint. In this section we will prove Theorem \ref{thm:main3}, which claims that $\chi(H_{\mathcal{F}})=\left\lceil{ \frac {\omega(H_{\mathcal{F}}) }{2}}\right\rceil$.

Note that a sub-family of intervals $\mathcal{L} \subset \mathcal{F}$ corresponds to a clique of $H_{\mathcal{F}}$ if and only if $\mathcal{L}$ is a pairwise disjoint subset of intervals of $\mathcal{F}$. Thus, $\omega(H_{\mathcal{F}})$ is the cardinality of the largest pairwise disjoint subset of intervals in $\mathcal{F}$. 

\textbf{Proof of Theorem \ref{thm:main3}}:
As we already know $\chi(H_{\mathcal{F}})\geq \left\lceil{ \frac {\omega(H_{\mathcal{F}}) }{2}}\right\rceil$ then it suffices to show that $\chi(H_{\mathcal{F}})\leq \left\lceil{ \frac {\omega(H_{\mathcal{F}}) }{2}}\right\rceil$. We proceed by induction 
on the size of  $\omega(H_{\mathcal{F}})$. Let $\omega(H_{\mathcal{F}}) = 2$ then the $3$-hypegraph $H_{\mathcal{F}}$ has an empty set of edges, thus $\chi(H_{\mathcal{F}})=1$.

For the remainder of the induction argument we will deal with the odd case and the even case separately.

Assume that $\omega(H_{\mathcal{F}}) = 2n + 1$. Let $\{L_0, L_1, ..., L_{2n} \} \subset \mathcal{F} $  be a pairwise disjoint set of intervals. Without loss of generality, we might assume that there is no other interval of ${\mathcal{F}}$ contained in $L_i$ for every $i=0,\dots 2n$. Denote by $G_0$ the set of intervals of ${\mathcal{F}}$ that intersect $L_0$. Note that, as no interval of ${\mathcal{F}}$ is contained in $L_0$, necessarily out of every three intervals in $G_0$ two of them intersect. Thus, $G_0$ corresponds to a independent set of $H_{\mathcal{F}}$. 

Next observe that, ${\mathcal{F}}  \setminus G_0$  (consisting of all intervals in $S^1$ that do not intersect  $L_0$), may be regarded as a set of intervals in the real line, assume with out lost of generality that they are ordered form left to right as $L_1,L_2,\dots ,L_{2n}$. As $\omega(H_{\mathcal{F}}) = 2n + 1$, the maximum number of pairwise disjoint intervals of ${\mathcal{F}}  \setminus G_0$ is $2n$. Consider $l_1,l_2,\dots ,\l_n$ the left hand side of each of this intervals. For $1 \leq i \leq n$, let $G_i$ be the subset of intervals of ${\mathcal{F}} \setminus G_0$ that intersects $\{l_{2i-1}, l_{2i}\}$ for each $i = 1, ...n$ . Clearly, each $G_i$ is an independent set because it satisfies that out of every three intervals in $G_i$, two intersect. Consequently, $H_{\mathcal{F}}$ contains $n+1$ independent sets and thus, $\chi(H_{\mathcal{F}}) \leq n + 1$.

Now assume that $\omega(H_{\mathcal{F}}) = 2n$. Let $\{L_1,L_2,...,L_{2n}\} \subset {\mathcal{F}}$ be a pairwise disjoint set of intervals. Without loss of generality, suppose that the distance between 
$L_{2n-1}$ and $L_{2n-2}$ is the smallest possible between all choices of pairwise disjoint subsets of ${\mathcal{F}}$. This implies two facts; firstly, that the collection, $G_0$ of intervals 
of ${\mathcal{F}}$ which intersect $\{L_{2n-1},L_{2n-2}\}$ satisfies that out of every three intervals in $G_0$, two intersect, thus $G_0$ is an independent set of the hypergraph $H_{\mathcal{F}}$; secondly that the collection of intervals ${\mathcal{F}}  \setminus G_0$ may be regarded as a set of intervals in the real line. Here we also have that the maximum number of pairwise 
disjoint intervals of ${\mathcal{F}}  \setminus G_0$ is $2n-2$, and as before assume that they are ordered from left to right as $L_1,L_2,\dots L_{2n-2}$ and each one of them do not contain completely any other interval, consider $l_1,l_2,\dots ,\l_{2n-1}$ the left hand side of each one of this intervals. Let $G_i$ be the subset of intervals of ${\mathcal{F}}  \setminus G_0$ that intersects $\{l_{2i-1}, l_{2i}\}$ for each $i = 1, ...n-1$.  Clearly, each $G_i$ is independent because it satisfies that out of every three intervals in $G_i$, two intersect. Consequently  $\chi(H_{\mathcal{F}}) \leq n $. This, concludes the proof of the theorem.
\qed

\section{Conclusions}

Clearly any sub-hypergraph of a circle interval $3$--hypergraph is a circle interval $3$--hypergraph, thus Theorem \ref{thm:main3} implies that any circle interval $3$--hypergraph, as well as all its induced sub-hypergraphs satisfy  $\chi(H)=\left\lceil \frac{n}{2}\right\rceil$ which, in the classic context of perfection would mean that this hypergraphs are ``perfect". 

Initially  we were aiming to find  ``perfection" on  comparability $3$--hypergraphs as a natural generalization of comparability graphs, fact that turns out to be false do to Theorem \ref{thm:main1}  and \ref{thm:main2}, however it remains an interesting  question what does perfection in hypergraphs means.

\medskip

\noindent 
{\bf Acknowledgments} 

N.  Garc\'{i}a-Col\'{i}n was partially supported by CONACyT: 166951, and
A. Montejano and  D. Oliveros partially supported by CONACyT: 166306 and PAPIIT: IN101912

Finally, we acknowledge support for Center of Innovation in Mathematics, CINNMA.





\bibliographystyle{model1a-num-names}



\end{document}